\documentclass{psapm-l}

\newtheorem{theorem}{Theorem}[section]
\newtheorem{lemma}[theorem]{Lemma}
\theoremstyle{definition}
\newtheorem{definition}[theorem]{Definition}
\newtheorem{example}[theorem]{Example}
\newtheorem{xca}[theorem]{Exercise}

\theoremstyle{remark}
\newtheorem{remark}[theorem]{Remark}

\numberwithin{equation}{section}

\newcommand{\abs}[1]{\lvert#1\rvert}

\newcommand{\blankbox}[2]{%
  \parbox{\columnwidth}{\centering
    \setlength{\fboxsep}{0pt}%
    \fbox{\raisebox{0pt}[#2]{\hspace{#1}}}%
  }%
}

\begin{document}

\title[Blow up of solutions to
the compressible magnetohydrodynamic]{Blow up of smooth solutions to
the barotropic compressible magnetohydrodynamic equations with
finite mass and energy}

\author{ Olga Rozanova}
\address{Department of Differential Equations \& Mechanics and Mathematics Faculty,
Moscow State University, Moscow, 119992,
 Russia}
\email{rozanova@mech.math.msu.su}
\thanks{Supported by  DFG 436 RUS
113/823/0-1 and \,RFBR 06-01-00441-a}


\date{January 1, 1994 and, in revised form, June 22, 1994.}



\begin{abstract}
We prove that the smooth solutions to the Cauchy problem for the
three-dimensional compressible barotropic magnetohydrodynamic
equations with  conserved total mass and finite total energy lose
the initial smoothness within a finite time. Further, we show that
the same result holds for the solution to the Cauchy problem for the
multidimensional compressible Navier-Stokes system. Moreover, for
the solution with a finite momentum of inertia we get the two-sided
estimates of different components of total energy.
\end{abstract}

\maketitle

\section{The finite time blow up result}
The set of equations which describe compressible viscous
magnetohydrodynamics are a combination of the compressible
Navier-Stokes equations of fluid dynamics and Maxwell's equations of
electromagnetism. We consider the system of partial differential
equations for the three-dimensional viscous compressible
magnetohydrodynamic flows in the Eulerian coordinates \cite{LL} for
barotropic case:
$$\partial_t \rho+{\rm div}_x (\rho u)=0,\eqno(1.1)$$
$$\partial_t(\rho u)+{\rm Div}_x (\rho u \otimes u)+\nabla_x p(\rho)=({\rm curl}_x H)\times H+{\rm Div_x} T,\eqno(1.2)$$
$$\partial_t H-{\rm curl}_x (u\times H)=-{\rm curl}_x(\nu\,{\rm curl}_x H),\qquad {\rm div}_x H=0, \eqno(1.3)     $$ where
$\rho, \, u=(u_1,u_2,u_3),\, p, \, H=(H_1,H_2,H_3),$ denote the
density, velocity, pressure and magnetic field. We denote ${\rm
Div}$ and $\rm div$ the divergency of tensor and vector,
respectively. Here $\,T\,$ is the stress tensor given by the Newton
law
$$T=T_{ij}=\mu \,(\partial_iu_j+\partial_j u_i)+\lambda \,{\rm div} u\, \delta_{ij},\eqno (1.4)$$
where the constants $\mu$ and $\lambda$ are the coefficient of
viscosity and the second coefficient of viscosity, $\nu\ge 0$ is the
coefficient of diffusion of the magnetic field.  We assume that $\mu
> 0, \,\lambda+\frac{2}{3}\mu> 0.$

The state equation has the form
$$
p=A\rho^\gamma.\eqno(1.5)
$$
Here $A>0,\,\gamma>1$ are constants.

Although the electric field $E$ does not appear in the MHD system
(1.1)-(1.3), it is indeed induced according to the  relation
$$E=\nu\,{\rm curl}_xH- u\times H$$
by the moving conductive flow in the magnetic field.

The MGD system (1.1) -- (1.5) is supplemented with the initial data
$$(\rho,u,H)\Bigl|_{t=0}=(\rho_0(x),u_0(x),H_0(x))\in{ \mathcal H}^m({\mathbb
R}^3),\,m>3.\eqno(1.6)$$


The local in time existence to the above Cauchy problem for the
density separated from zero follows from \cite{VolpertKhudiaev},
however the results on global solvability are very scanty. For
example, even for the one-dimensional case, the global existence of
classical solutions to the full perfect MHD equations with large
data remains unsolved (the smooth global solution near the constant
state in one-dimensional case is investigated in \cite{KO}). The
existence of global weak solutions was proved recently \cite{HW2}
(see also \cite{HW1} for references).


\vskip1cm

The system (1.1)-- (1.3) is the differential form of balance laws
for the material volume $\Omega(t)$; it expresses  conservation of
mass
$$m=\int\limits_{\Omega(t)}\rho \, d x,$$  balance of momentum
$$P=\int\limits_{\Omega(t)}\rho u\, d x,$$ and balance of total energy
$$\mathcal E=\int\limits_{\Omega(t)}\left(\frac{1}{2}\rho
|u|^2+\frac{|H|^2}{2}+\frac{A\rho^\gamma}{\gamma-1}\right)\, d x
\,=E_{k}(t)+E_m(t)+E_{i}(t).$$ Here $E_k(t),\,$ $E_m(t)$ and
$E_i(t)$ are the kinetic, magnetic and internal components of
energy, respectively.

 If we regard $\Omega (t)={\mathbb R}^3,$ the conservation of
mass, momentum and non-increasing of total energy takes place
provided the components of solution  decrease at infinity
sufficiently quickly.

\vskip1cm

Our main result is the following:

\begin{theorem} Assume $\gamma\ge\frac{6}{5}.$ If  the momentum $P\ne 0,$ then there exists no global in time
classical solution to (1.1) -- (1.5), (1.6) with conserved mass and
finite total energy.
\end{theorem}

If $H\equiv 0,$ we get the barotropic system of the compressible
Navier-Stokes equations. We can consider it in any space dimensions
without physical restriction $n=3.$ In this case the viscosity
coefficients should satisfy the inequalities $\mu
> 0, \,\lambda+\frac{2}{n}\mu> 0.$ The following theorem is an
extension of the result of \cite{RozJDE} to the larger class of
solutions.

\begin{theorem} Assume $n\ge 3$ and $\gamma\ge\frac{2n}{n+2}.$ If  the momentum $P\ne 0,$ then there exists no global
in time classical solution to the barotropic system of the
compressible Navier-Stokes equations with conserved mass and finite
total energy.
\end{theorem}

The next lemma is proved in \cite{RozJDE} for both barotropic and
non-barotropic cases.

\begin{lemma} Let $n\ge 3,\,\gamma\ge \frac{2n}{n+2}$ and $u\in {\mathcal H}^1({\mathbb R}^n).$ If $|P|\ne 0,$ then there
exists a positive constant $K$ such that for the solutions with
conserved mass and finite total energy the following inequality
holds:
$$
\int\limits_{{\mathbb R}^n}\,|Du|^2\,dx\ge
\,K\,E_i^{-\frac{n-2}{n(\gamma-1)}}(t). \eqno(1.7)
$$
\end{lemma}

In the dimension $n=3$ it has the form:

\begin{lemma} Let $\gamma\ge \frac{6}{5},$\,and $u\in {\mathcal H}^1({\mathbb R}^n).$ If $|P|\ne 0,$ then there
exists a positive constant $K$ such that for classical solutions to
(1.1)-- (1.5) with conserved mass and total finite energy the
following inequality holds:
$$
\int\limits_{{\mathbb R}^3}\,|Du|^2\,dx\ge
\,K\,E_i^{-\frac{1}{3(\gamma-1)}}(t). \eqno(1.8)
$$
\end{lemma}

{\it Proof of Lemma 1.4} First of all we use the H\"older inequality
to get
$$
|P|=\Big|\int\limits_{{\mathbb R}^3}\,\rho u \,dx\Big|\le
\left(\int\limits_{{\mathbb R}^3}\,\rho^{\frac{6}{5}}
\,dx\right)^{\frac{5}{6}}\left(\int\limits_{{\mathbb R}^3}\, |u|^{6}
\,dx\right)^{\frac{1}{6}}.\eqno(1.9)
$$
Further, using the Jensen inequality we have for $5(\gamma-1)\ge 1$
(or $\gamma\ge\frac{6}{5}$)
$$
\left(\frac{1}{m}\,\int\limits_{{\mathbb R}^3}\,\rho^{\frac{6}{5}}
\,dx\right)^{5(\gamma-1)}\le \frac{\int\limits_{{\mathbb
R}^3}\rho^\gamma \,dx}{m}=\frac{(\gamma-1) E_i(t)}{mA}.
$$
Thus,
 the latter inequality and (1.9) give
$$
|P|\le K_1\, \left(E_i(t)\right)^{\frac{1}{6(\gamma-1)}}
\left(\int\limits_{{\mathbb R}^3}\, |u|^{6}
\,dx\right)^{\frac{1}{6}},\eqno(1.10)
$$
with the positive constant $K_1$ that depends on $\gamma$ and $ m.$
Further, we take into account the inequality
$$
\left(\int\limits_{{\mathbb R}^n}\,
|u|^{\frac{2n}{n-2}}\,dx\right)^{\frac{n-2}{n}} \le K_2
\,\int\limits_{{\mathbb R}^n}\, |Du|^2 \,dx,\eqno(1.11)
$$
where the constant $K_2>0 $ depends on $n,$  $\,n\ge 3.$ The latter
inequality holds for $u\in H^1({\mathbb R}^n)$ (\cite{Hebey}, p.22)
and follows from the Sobolev embedding.

Thus, the lemma statement follows from (1.10), (1.11), with the
constant $K=\frac{|P|^2 }{K_1^2 K_2}.$

\vskip1cm

{\it Proof of Theorem 1.1.}
Let us note that the total energy $\mathcal E$ satisfies the
inequality
$$\mathcal E'(t)\le -\nu \int \limits_{{\mathbb R}^3}|{\rm curl}_x H|^2\, dx-
\sigma \int \limits_{{\mathbb R}^3}|D u|^2\, dx\le 0,\eqno(1.12)$$
with some positive constant $\sigma,$ therefore $\mathcal E(t)\le
\mathcal E(0).$

From (1.8) we have
$$\mathcal E'(t)\le -\sigma \,K \,(E_i(t))^{-\frac{1}{3(\gamma-1)}}\le -\sigma \,K \,
({\mathcal E}(0))^{-\frac{1}{3(\gamma-1)}}=const<0,\eqno (1.13)  $$
this contradicts to the non-negativity of $\mathcal E(t).$ Thus,
Theorem 1.1 is proved.

{\it Proof of Theorem 1.2} is similar. Estimate (1.7) holds in any
dimensions greater or equal than 3, therefore we can apply Lemma 1.3
to the case $H\equiv 0$ to obtain
$$\mathcal E'(t)\le -\sigma \,K \,(E_i(t))^{-\frac{n-2}{n(\gamma-1)}}\le
-\sigma \,K \,({\mathcal E}(0))^{-\frac{n-2}{n(\gamma-1)}}=const<0.
$$

\section{Estimates of energy for solutions with finite moment of inertia}

In this section we are going to obtain more specified estimate for
the decay rate of total energy for more restricted classes of
solutions.

 Let us introduce the  functionals
$$G(t)=\frac{1}{2} \int\limits_{{\mathbb
R}^3}\rho(t,x)|{ x}|^2\,dx,\qquad F(t)=\int\limits_{{\mathbb R}^3}
({u},{ x})\rho\,dx,$$ where the first one is the momentum of
inertia, the inner product of vectors is denoted as $(.,.)$ .

Taking into account the general Stokes formula one can readily
calculate that for classical compactly supported solution to
(1.1)--(1.5) the following identities hold:
$$G'(t)=F(t),\eqno(2.1)$$
$$F'(t)=2 E_k(t)+ E_m(t) + 3(\gamma-1) E_i(t).\eqno(2.2)$$
In fact, formulas are obtained in \cite{Chemin} for gas dynamics
equations with finite moment of inertia (not especially compactly
supported). However, evidently that  (2.1), (2.2) are satisfied for
solutions decaying at infinity sufficiently quickly together with
their derivatives. One can impose a special rate of decay on the
solution and then prove (2.1), (2.2) as a lemma (as in
\cite{RozJDE}). The alternative way is to take (2.1), (2.2) as a
definition of the rate of decay.

\begin{definition}
We will call a solution $(\rho, u,H)$ to the Cauchy problem (1.1) --
(1.5), (1.6) {\it highly decreasing at infinity}, if $\rho(t,x)\ge
0$ and the identities (2.1), (2.2) hold.
\end{definition}

Evidently, the integrals of mass, momentum and total energy are
finite for highly decreasing at infinity solutions at any fixed
moment of time.

\medskip

We firstly get two-sided estimates of $G(t).$

\begin{lemma} For classical highly decreasing at infinity solution to (1.1) -- (1.5)
the following estimates hold:

if $\gamma\le \frac{5}{3}, $ then
$$\frac{P^2}{2m}t^2+F(0)t+G(0)\le G(t)\le  \mathcal E(0) t^2+F(0)t+G(0);\eqno(2.3)$$
if $\gamma> \frac{5}{3}, $ then
$$\frac{P^2}{2m}t^2+F(0)t+G(0)\le G(t)\le \frac{3(\gamma-1)}{2} \mathcal E(0) t^2+F(0)t+G(0).\eqno(2.4)$$
\end{lemma}

{\it Proof of Lemma 2.2} Let us note that  (1.12) implies $\mathcal
E(t)\le \mathcal E(0).$

 Identities (2.1) and (2.2) imply for $\gamma\le\frac{5}{3}$
$$
G''(t)=2E_k(t)+E_m(t)+3(\gamma-1)E_i(t)=
2E_k(t)+E_m(t)+3(\gamma-1)E_i(t)=
$$
$$=2\mathcal E (t) - E_m(t)-
(2-3(\gamma-1))E_i(t)\le 2\mathcal E (0),\eqno(2.5)
$$
and for $\gamma>\frac{5}{3}$
$$
G''(t)=2E_k(t)+E_m(t)+3(\gamma-1)E_i(t)=
2E_k(t)+E_m(t)+3(\gamma-1)E_i(t)=
$$
$$=3(\gamma-1)\mathcal E (t) -(3\gamma-4) E_m(t)-
(3\gamma-5)E_k(t)\le 3(\gamma-1)\mathcal E (0).\eqno(2.6)
$$
Inequality (2.5) and (2.6) result in the right-hand sides in (2.3)
and (2.4).

 Further, since $E_k(t)\ge\frac{P^2}{2m}$ (that follows
from the H\"older inequality),  $E_i(t)$ and $E_m(t)$ are
nonnegative, then from (2.1), (2.2) we obtain
$$G''(t)\ge \frac{P^2}{m}.\eqno(2.7)$$ Thus, the left-hand sides estimates in (2.5) and (2.6) are obtained.
Lemma 2.2 is proved.


The next step is two-sides estimate of $E_i(t)+E_m(t).$

\begin{lemma} If $P\ne 0,$ then
 for highly decreasing at infinity classical solutions to (1.1) -- (1.5)
for sufficiently large $t$ the following estimates hold:
$$\frac{C_1}{G^{3(\gamma-1)/2}(t)}\le E_i(t) +E_m(t)\le \frac{C_2}{G^{3(\gamma-1)/2}(t)},\eqno(2.8)$$
for $\gamma\le \frac{4}{3},$ and
$$\frac{C_1}{G^{3(\gamma-1)/2}}\le E_i(t)+E_m(t) \le \frac{C_2}{G^{1/2}(t)},\eqno(2.9)$$
for $\gamma > \frac{4}{3},$ with positive constants $C_1, C_2
\,(C_1\le C_2).$
\end{lemma}

{\it Proof of Lemma 2.3} The lower estimate is due to \cite{Chemin}.
For an arbitrary dimension of space it follows from the inequality
$$
\|f\|_{L^1({\mathbb R}^n;\,dx)}\,\le\,C_{\gamma,
n}\|f\|^{\frac{2\gamma}{(n+2)\gamma-n}}_{L^\gamma({\mathbb
R}^n;\,dx)}\,\|f\|^{\frac{n(\gamma-1)}{(n+2)\gamma-n}}_{L^1({\mathbb
R}^n;\,|x|^2\,dx)},
$$
with
 $$C_{\gamma,
n}=\left(\frac{2\gamma}{n(\gamma-1)}\right)^
{\frac{n(\gamma-1)}{(n+2)\gamma-n}} +
\left(\frac{2\gamma}{n(\gamma-1)}\right)^
{\frac{-2\gamma}{(n+2)\gamma-n}}.$$

Thus the constant
$$C_1=\frac{A}{\gamma-1}(m C_{\gamma,
n}^{-1})^\frac{\gamma(n+2)-n}{2}.$$

The method of the upper estimate of $E_i(t)+E_m(t)$ is also similar
to \cite{Chemin}. Namely, let us consider the function
$Q(t)=4G(t)\mathcal E(t)-F^2(t).$ The H\"older inequality gives
$F^2\le 4 G(t) E_k(t),$ therefore $\mathcal
E(t)=E_k(t)+E_m(t)+E_i(t)\ge E_i(t)+\frac{F^2(t)}{4G(t)}$ and
$$
E_i(t)+E_m(t)\le \frac {Q(t)}{4 G(t)}.\eqno(2.10)
$$
We notice also that $Q(t)>0$ provided the pressure does not equal to
zero identically. Then taking into account (2.1), (2.2) and (1.12)
we have
$$
Q'(t)=4 G'(t) \mathcal E(t)-2 G'(t) G''(t)+ 4G(t)\mathcal
E'(t)=$$$$= 2\,G'(t)\,((5-3\gamma)\,E_i(t)+E_m(t))+4G(t)\mathcal
E'(t)\le $$$$\le 2\,G'(t)\,((5-3\gamma)\,E_i(t)+E_m(t)).\eqno(2.11)
$$

Since $P\ne 0,$ inequality (2.7)  guarantees  the positivity of
$G'(t)$ for sufficiently large $t,$ then for $\gamma\le\frac{4}{3}$
$$Q'(t)=
2\,G'(t)\,((5-3\gamma)\,(E_i(t)+E_m(t))+(3\gamma-4)E_m(t))\le$$
$$\le
2\,G'(t)\,(5-3\gamma)\,(E_i(t)+E_m(t));\eqno(2.12)
$$
and for $\gamma>\frac{4}{3}$
$$Q'(t)=
2\,G'(t)\,((E_i(t)+E_m(t))+(4-3\gamma)E_m(t))\le
2\,G'(t)\,(E_i(t)+E_m(t)).\eqno(2.13)
$$
Thus, for $\gamma\le \frac{4}{3}$ inequality (2.12) implies
$$\frac{Q'(t)}{Q(t)}\le \frac{5-3\gamma}{2}\,\frac{G'(t)}{G(t)},\eqno(2.14)$$
whereas for $\gamma>\frac{4}{3}$ inequality (2.13) results in
$$\frac{Q'(t)}{Q(t)}\le \frac{1}{2}\,\frac{G'(t)}{G(t)}.\eqno(2.15)$$

Therefore (2.10), (2.14) and (2.15) imply the upper estimates in
(2.8), (2.9), where
$C_2=\frac{1}{4}\,Q(0)\,G^{\frac{3\gamma-5}{2}}(0)$ in (2.8) and
$C_2=\frac{Q(0)}{4\sqrt{G(0)}}$ in (2.9). The proof of Lemma 2.3 is
over.

\vskip1cm

Now we can estimate the rate of decay of total energy $\mathcal E.$
From (1.12), (1.8) we have
$$\mathcal E'(t)\le -\sigma \,K \,(E_i(t))^{-\frac{1}{3(\gamma-1)}}.\eqno (2.16)  $$
Further, nonnegativity of $E_m(t),$  Lemmas 2.2 and 2.3 and
inequality (2.16) imply
$$
\mathcal E'(t)\le-\nu\, K\,
C_2^{-\frac{1}{3(\gamma-1)}}\,G^{\frac{1}{6}}(t)\le - L
\,t^{\frac{1}{3}},
$$
with a positive constant $L,$ for $\gamma\le \frac{4}{3},$ and
$$
\mathcal E'(t)\le-\nu\, K\,
C_2^{-\frac{1}{3(\gamma-1)}}\,G^{\frac{1}{6(\gamma-1)}}(t)\le - L \,
t^{\frac{2}{3(\gamma-1)}},
$$
for $\gamma> \frac{4}{3}.$ In both cases this contradicts to the
non-negativity of $\mathcal E(t).$ However, the decay of total
energy in the case of highly decreasing solutions is greater that
follows from inequality (1.13).

\vskip 0.5cm

The author thanks Prof. Ronghua Pan for attraction of attention to
this problem and helpful discussion.

\end{document}